\documentclass[12pt]{article} 
\usepackage{amssymb,latexsym}
\def\neweq{\setcounter{equation}{0}}
\newtheorem{theorem}{Theorem}[section]
\newtheorem{pr}[theorem]{Proposition}
\newtheorem{cor}[theorem]{Corollary}
\newtheorem{de}[theorem]{Definition}

\newtheorem{lem}[theorem]{Lemma}

\font\germ=eufm10
\def\F{\hbox{\germ F}}

\def\S{\mathfrak{S}}
\def\Q{\mathbb{Q}}

\def\s{\hbox{\germ S}}
\def\<{\langle}
\def\>{\rangle}

\def\N{\mathbb{Z}}
\def\S{\mathbb{S}}
\def\Z{\mathbb{Z}}

\def\wt{\widetilde}

\def\ds{\displaystyle}
\def\qed{\hfill$\vrule height 2.5mm width 2.5mm depth 0mm$}

\title{ $t$--Deformations of quantum Schubert polynomials}
\date{ }
\author{Anatol N. Kirillov  \\~ \\
{\small {\it CRM, University of Montreal}} \\
{\small {\it C.P. 6128, Succursale A, Montreal (Quebec), H3C 3J7, 
Canada}}\\
{\small {\it and}}
 \\
{\small {\it Steklov Mathematical Institute,}} \\
{\small {\it Fontanka 27, St.Petersburg, 191011, Russia}}}

\begin{document}
\maketitle
\hskip 9cm{\it For my daughter's 14-th birthday}

\begin{abstract}
We construct a certain solution ${\cal F}(t)$ to the 
Witten--Dijkgraf--Verlinde--Verlinde equation 
related to the small quantum cohomology ring of flag variety, and study the
$t$--deformation of quantum Schubert polynomials corresponding to 
this solution.
\end{abstract}

\section{Introduction}
\label{sec:intro}
\neweq

The cohomology ring of the flag variety $Fl_n=SL_n/B$ is isomorphic to 
the quotient of the polynomial ring $P_n:=\Z [x_1,\ldots ,x_n]$ by 
the ideal $I_n$ generated by symmetric polynomials without constant term:
$$H^*(Fl_n,\Z )\cong P_n/I_n.
$$
The Schubert cycles $X_w:=\overline{BwB/B}$, $w\in W$, form a linear basis of 
the homology group $H_*(Fl_n,\Z )$ and via the Poincare duality they are 
represented by the (geometric) Schubert polynomials $X_w(x)$ in the 
cohomology ring $H^*(Fl_n,\Z )$. It was discovered by A.~Lascoux and 
M.-P.~Sch\"utzenberger, [LS1], [LS2], that there exists the set of
distinguish representatives $\s_w(x)\in P_n$ of the geometric Schubert 
polynomials $X_w(x)\in P_n/I_n$ with nice algebraic, combinatorial and 
geometric properties. Follow A.~Lascoux and \hbox{M.-P.~Sch\"utzenberger} 
[LS1] 
these distinguish representatives $\s_w(x)$ are called {\it Schubert
polynomials}. These polynomials form a stable, homogeneous, orthonormal 
basis in the ring of polynomials $P_n=\Z 
[x_1,\ldots ,x_n]$ indexed by permutations $w\in S^{(n)}=\{ w\in 
S_{\infty},~l(c(w))\le n\}$. We refer the reader to [M] and [F] for 
detailed account on Schubert polynomials. One of the most deep and 
fundamental properties of the Schubert polynomials $\s_w(x)$ is that 
their structural constants $c_{uv}^w$, defined from the decomposition
$$\s_u\s_v=\sum_wc_{uv}^w\s_w~{\rm mod}~I_n,
$$
appear to be nonnegative integers. The only known proof of this fact 
appeals to algebraic geometry and is based on an interpretation of the 
structural constant $c_{uv}^w$ as the intersection number of the Schubert 
cycles $X_u,X_v$ and $X_{w_0w}$. One of the most fundamental 
problems of the Schubert Calculus is to find an algebraic proof of 
nonnegativity of the structural constants $c_{uv}^w$ and give their 
combinatorial interpretation (Littlewood--Richardson's problem for 
Schubert polynomials).

Since the {\it quantum cohomology} of compact 
complex Kahler manifolds, was introduced by C.~Vafa [V], the computation 
of the (big) quantum cohomology 
rings of flag varieties became a priority problem of {\it Quantum Schubert 
Calculus}. We refer the reader to [FP] and [MS] where the definition and 
basic properties of quantum cohomology may be found. 

The essential new ingredient of the quantum cohomology theory in 
comparison with the classical one is the existence of the Gromov--Witten 
potential ${\cal F}(t)$ which contains all information about the 
multiplication rules in the quantum cohomology ring. The Gromov--Witten 
potential ${\cal F}(t)$ satisfies the celebrated 
Witten--Dijkgraaf--Verlinde--Verlinde equation (WDVV--equation):
\begin{equation}\label{0.1}
\sum_{\nu ,\mu}\frac{\partial^3{\cal F}}{\partial t_i\partial t_j\partial 
t_{\nu}}g^{\nu\mu}\frac{\partial^3{\cal F}}{\partial t_{\mu}\partial 
t_k\partial t_l}=\sum_{\nu ,\mu}\frac{\partial^3{\cal F}}{\partial 
t_i\partial t_l\partial t_{\nu}}g^{\nu\mu}\frac{\partial^3{\cal 
F}}{\partial t_{\mu}\partial t_j\partial t_k},
\end{equation}
for all $i,j,k,l$, where $g^{\nu\mu}$ stands for the inverse 
matrix of $g_{\mu\nu}$, where $g_{\mu\nu}$ is the intersection matrix:
$$g_{\nu\mu}=\ds\frac{\partial^3{\cal F}}{\partial t_0\partial t_{\nu}
\partial t_{\mu}}.
$$

In particular case of the flag variety $Fl_n$, the WDVV--equation takes 
the following form
\begin{equation}\label{0.2}
\sum_{v\in S_n}\frac{\partial^3{\cal F}}{\partial t_{w_1}\partial t_{w_2}
\partial t_v}\cdot\frac{\partial^3{\cal F}}{\partial t_{w_0v}
\partial t_{w_3}\partial t_{4}}=\sum_{v\in S_n}\frac{\partial^3{\cal F}}
{\partial t_{w_2}\partial t_{w_3}\partial t_v}\cdot\frac{\partial^3{\cal 
F}}{\partial t_{w_0v}\partial t_{w_1}\partial t_{w_4}},
\end{equation}
for all $w_1,w_2,w_3,w_4\in S_n$, where $w_0$ is the longest element of
$S_n$, and $t=(t_w,~w\in S_n)$ stands for the set of 
independent variables $t_w$ parameterized by permutations $w\in S_n$. 

It was conjectured in [KM] that the WDVV--equation (\ref{0.2}) 
defines uniquely the Gro\-mov--Witten potential ${\cal F}(t)$ for the flag 
variety $Fl_n$ if the following conditions are satisfied:

1) Normalization:
$$\frac{\partial^3{\cal F}}{\partial t_1\partial t_v\partial 
t_w}=\delta_{v,w_0w};
$$

2) Initial conditions:
$$\frac{\partial^3{\cal F}}{\partial t_{s_k}\partial t_{s_k}\partial 
t_{w_0}}=q_k,~~~1\le k\le n-1,
$$
where for each $k\in [1,n-1]$, $q_k$ is a certain constant, and 
$s_k=(k,k+1)$ denotes the transposition that 
interchanges $k$ and $k+1$, and fixes all other elements of $[1,n]$;

3) Degree conditions:
$$\frac{\partial^3{\cal F}}{\partial t_u\partial t_v\partial t_w}=0,
$$
if either $l(u)+l(v)+l(w)<l(w_0)$ or difference $l(u)+l(v)-l(w_0)$ is an 
odd positive integer.

In particular case when $t_w=0$ for all $w\in S_n$ such that $l(w)\ge 2$, 
the quantum cohomology ring ${\rm QH}^*(Fl_n)$ (the so--called small 
quantum cohomology ring)  was computed by A.~Givental and B.~Kim [GK], 
see also [C1]. The study of (small) quantum Schubert polynomials 
corresponding to the small quantum cohomology ring ${\rm QH}^*(Fl_n)$  of 
the flag variety $Fl_n$ was initiated in [FGP] and independently in [KM], 
see also [C2]. The Grassmannian case was considered earlier in [B], [C1] 
and [W].

The main goal of the present paper is to construct a toy model for 
big quantum cohomology ring and big quantum Schubert polynomials for the 
flag variety $Fl_n$. More precisely, to construct the (toy) Gromov--Witten
potential ${\cal F}(t)$ which appears to be a solution to the WDVV--equation
(\ref{0.2}). Based on this solution, we define the $t$--deformation 
$\wt\s_w^t$ of the quantum Schubert polynomials $\wt\s_w$ and investigate 
some of their properties. The (toy) Gromov--Witten potential ${\cal F}(t)$ 
is defined to be the 
Grothendieck residue of the function $U(x)=\exp\left(\ds\sum_{w\in 
S_n}t_w\wt\s_w(x)\right)$ with respect to the ideal $\wt I_n$:
$${\cal F}(t)=\langle\exp\left(\sum_{w\in 
S_n}t_w\wt\s_w(x)\right)\rangle_{\wt I_n}.
$$
Here the function ${\cal F}(t)$ is a natural generalization of the quantum 
generating volume function $V(z;q)$ introduced by A.~Givental and B.~Kim 
[GK].

The content of the paper is arranged as follows. In Section~2 we review 
the definition and some basic properties of the quantum Schubert 
polynomials to be used in Sections~3 and 4. In Section~3 we 
introduce and study the (toy) Gromov--Witten potential ${\cal F}(t)$
and the $t$--deformation of quantum Schubert polynomials related to
the potential ${\cal F}(t)$. In Section~4 we construct a certain Lax pair 
related to yet another deformation $X_w^t$ of the quantum Schubert
polynomials.

\bigskip

\textsc{Acknowledgments.} I am thankful to Alex Kasman and Toshiaki Maeno 
for fruitful discussions.

\newpage
\section{Quantum Schubert polynomials}
\label{sec:qsp}
\neweq

In this section we review the definition and basic properties of the 
quantum Schubert polynomials. The study of (small) quantum Schubert 
polynomials for flag variety $Fl_n$ was initiated in [FGP] and, 
independently, in [KM]. The case of Grassmannian varieties was considered 
earlier in [B], [C1] and [W]. In our exposition we follow [KM].

Let $X_n=(x_1,\ldots ,x_n)$, $Y_n=(y_1,\ldots ,y_n)$ be two sets of 
independent variables, put
$${\wt\s}_{w_0}(x,y):={\wt\s}_{w_0}^{(q)}(X_n,Y_n)=\prod_{i=1}^{n-1}
\Delta_i(y_{n-i}~|~X_i),
$$
where $\Delta_k(t|X_k):=\ds\sum_{j=0}^kt^{k-j}e_j(X_k|q_1,\ldots 
,q_{k-1})$ is the generating functions for the quantum elementary 
symmetric polynomials in the variables $X_k=(x_1,\ldots ,x_k)$, i.e.

\vskip 0.3cm
$\Delta_k(t|X_k):=\sum_{i=0}^ke_i(X_k|q)t^{k-i}$
\vskip 0.2cm
\begin{equation}\label{1.1}
=\det\left(\begin{array}{ccccccc} x_1+t & q_1&0 &\ldots &\ldots 
&\ldots &0\\
-1 & x_2+t & q_2 &0 &\ldots &\ldots & 0\\
0 & -1 &x_3+t & q_3 & 0 &\ldots & 0 \\
\vdots &\ddots &\ddots & \ddots &\ddots &\ddots &\vdots \\
0&\ldots & 0 &-1&x_{k-2}+t &q_{k-2} & 0 \\
0 &\ldots &\ldots & 0 &-1 & x_{k-1}+t & q_{k-1}\\
0 &  \ldots &\ldots &\ldots & 0 & -1 & x_k+t
\end{array}\right) .
\end{equation}
\vskip 0.3cm
The determinant (\ref{1.1}) was introduced in connection to the (small) 
quantum cohomology ring of flag varieties by A.~Givental and B.~Kim [GK].

\begin{de}\label{d1.1} {\rm ([KM])} For each permutation $w\in S_n$, the 
quantum double Schubert polynomial ${\wt\s_w}(x,y)$ is defined to be
$${\wt\s}_w(x,y)=\partial_{ww_0}^{(y)}{\wt\s}_{w_0}(x,y),
$$
where divided difference operator $\partial_{ww_0}^{(y)}$ acts on the $y$ 
variables.
\end{de}
\begin{de}\label{d1.2} {\rm (cf. [KM])} For each permutation $w\in S_n$, 
the (small) quantum Schubert polynomial ${\wt\s}_w(x)$ is defined to be 
the $y=0$ specialization of the quantum double Schubert polynomials 
${\wt\s}_w(x,y)$:
$${\wt\s}_w(x,y)=\partial_{ww_0}^{(y)}{\wt\s}_{w_0}(x,y)|_{y=0}.
$$
\end{de}

Now we are going to review some basic properties of (small) quantum 
Schubert polynomial which will be used in the next sections. We start with 
reminding the result of A.~Givental and B.~Kim [GK], and 
I.~Ciocan-Fontanine [C1], on the structure of the small quantum 
cohomology ring ${\rm QH}^*(Fl_n)$ of the flag variety $Fl_n$:
$${\rm QH}^*(Fl_n)\cong\Z [x_1,\ldots ,x_n,q_1,\ldots ,q_{n-1}]/{\wt I_n},
$$
where the ideal $\wt I_n$ is generated by the quantum elementary symmetric 
polynomials\break $\wt e_i(x):=e_i(X_n|q_1,\ldots ,q_{n-1})$, $1\le i\le n$, 
with generating function $\Delta_n(t|X_n)$, see (\ref{1.1}). 

There exists a natural pairing $\langle f,g\rangle_Q$ on the ring of 
polynomials $\Z [X_n,q_1,\ldots ,q_{n-1}]$, and the small quantum 
cohomology ring ${\rm QH}^*(Fl_n)\cong\Z [X_n,q_1,\ldots ,q_{n-1}]
/{\wt I_n}$ which is 
induced by the Grothen\-dieck  residue with respect to the ideal 
${\wt I}_n$, see, e.g., [GK], [KM]:
$$\langle f,g\rangle_Q={\rm Res}_{\wt I_n}(f,g),~~~f,g\in\Z [X_n,q_1,\ldots 
,q_{n-1}].
$$
The pairing $\langle f,g\rangle_Q$ satisfies the following properties

1) ~$\langle f,g\rangle_Q=0$ if $f\in\wt I_n$;

2) ~$\langle f,g\rangle_Q$ defines a nondegenerative pairing in 
$\Z [x,q]/{\wt I}_n\cong{\rm QH}^*(Fl_n)$.

Let us denote by $\Lambda_n[q]=\Lambda_n\otimes_{\Z}\Z [q_1,\ldots 
,q_{n-1}]$, ${\cal H}_n[q]={\cal H}_n\otimes_{\Z}\Z [q_1,\ldots 
,q_{n-1}]$, where $\Lambda_n=\Z [x_1,\ldots ,x_n]^{\S_n}$ stands for the 
ring of symmetric functions, and ${\cal H}_n$ stands for the $\Z$--span 
generated by monomials $x^I=x_1^{i_1}\cdots x_n^{i_n}$, where 
$I\subset\delta_n=(n-1,n-2,\ldots ,1,0)$.

\begin{pr}\label{p1.1} {\rm ([FGP], [KM])}

$\bullet$ The quantum Schubert polynomials $\wt\s_w$, $w\in S_n$, form a 
$\Lambda_n[q]$--basis of $P_n\otimes_{\Z}\Z [q_1,\ldots ,q_{n-1}]$.

$\bullet$ The quantum Schubert polynomials $\wt\s_w$, $w\in S_n$, form a 
$\Z [q_1\ldots ,q_{n-1}]$--linear basis of ${\cal H}_n[q]$.

$\bullet$ (Orthogonality) Let $u,v\in S_n$. Then
$$\langle\wt\s_u,\wt\s_v\rangle_Q=\left\{\begin{array}{rl} 1, & {\rm 
if} \ \ u=w_0v,\\ 0, & {\rm otherwise.}\end{array}\right.
$$

$\bullet$ (Gromov--Witten's invariants of order 3) Let $u,v,w\in S_n$, 
then there exists a polynomial~ $\wt c_{uv}^w(q)\in \N [q_1,\ldots 
,q_{n-1}]$ with nonnegative integer coefficients such that
$$\wt\s_u\wt\s_v\equiv\sum_w\wt c_{uv}^w(q)\wt\s_w~{\rm mod}~\wt I_n.
$$
\end{pr}

Let us remind that the coefficient $\wt c_{uvd}^w$ of the polynomial
$$\wt c_{uv}^w(q)=\sum_{d\ge 0}\wt c_{uvd}^wq^d
$$
is called the Gromov--Witten invariant of order 3 and degree $d$ 
corresponding to the Schubert cycles $X_u$, $X_v$, and $X_w$.

\section{$t$--Deformation}
\label{sec:td}
\neweq

In this Section we introduce and study the (toy) Gromov--Witten potential 
${\cal F}(t)$ and define the $t$--deformation of the quantum Schubert 
polynomials corresponding to this potential.

To start, let us consider the function
\begin{equation}\label{2.1}
\exp\left(\sum_{w\in S_n}t_w\wt{\s}_w(x)\right)\in\Q [[X,q_1,\ldots 
,q_{n-1},\{ t_w\}_{w\in S_n}]],
\end{equation}
and denote by $K(x,t)=\ds\left[\exp\left(\sum_{w\in S_n}t_w\wt{\s}_w(x)
\right)\right]$ the projection of the function (\ref{2.1}) on the 
subspace $\overline{\cal H}_n(t,q):={\cal H}_n\otimes_{\Z}\Z [[q_1,\ldots
,q_{n-1}, \{ t_w\}_{w\in S_n}]]$ generated 
by the quantum Schubert polynomials $\wt{\s}_w(x)$, $w\in W$, which called
{\it kernel}. Thus, we have 
\begin{equation}\label{2.2}
K(x,t)=\sum_{w\in S_n}\varphi_w(t)\wt{\s}_{ww_0}(x),
\end{equation}
where coefficients $\varphi_w(t)$ belong to the formal 
power series ring $\Z [[q_1,\ldots ,q_{n-1},~\{ t_w\}_{w\in S_n}]]$.

\begin{de}\label{d2.1} Define the (toy) Gromov--Witten potential ${\cal 
F}(t)$ to be the Grothendieck residue of the kernel $K(x,t)$ with respect to 
the ideal $\wt I_n$: \ \ ${\cal F}(t)=\langle K(x,t)\rangle_{\wt I_n}$.
\end{de}

In other words, ${\cal F}(t)=\varphi_{\rm id}(t)$, i.e.
$$K(x,t)={\cal F}(t)\wt{\s}_{w_0}(x)+\sum_{w>1}\varphi_{w}(t)\wt{\s}_{ww_0}(x).
$$

\begin{lem}\label{l2.1} For each $w\in S_n$,
\begin{equation}\label{2.3}
\frac{\partial}{\partial t_w}{\cal F}(t)=\varphi_w(t).
\end{equation}
\end{lem}

{\it Proof.}~~ By definition,
\begin{eqnarray*}
\frac{\partial}{\partial t_w}{\cal F}(t)&=&\langle\frac{\partial}{\partial t_w}
\exp\left(\sum_{w\in S_n}t_w\wt{\s}_w(x)\right)\rangle_{\wt I}=\\
&=&\langle\wt{\s}_w(x)\cdot K(x,t)\rangle_{\wt I}=
\langle K(x,t),\wt{\s}_w(x)\rangle_Q=\varphi_w(t).
\end{eqnarray*}
\qed

It follows from Lemma~\ref{l2.1}, that
$$K(x,t)=\left(\sum_{w\in S_n}\frac{\partial}{\partial t_w}
\wt{\s}_{ww_0}(x)\right)\cdot {\cal F}(t).
$$
From now we will use the notation $[f]$, 
$f\in\Z[X_n][[q_1,\ldots ,q_{n-1},~\{ t_w\}_{w\in S_n}]]$ to denote the
projection of $f$ on the subspace ${\cal H}_n(t,q):={\cal H}_n\otimes_{\Z}
\Z [[q_1,\ldots ,q_{n-1},~\{ t_w\}_{w\in S_n}]]$ generated by the small 
quantum Schubert polynomials $\wt\s_w$, $w\in S_n$.

\begin{de}\label{d2.2} Define the $t$--deformation $\wt{\s}_w^t(x)$ of 
the quantum Schubert polynomial $\wt{\s}_w(x)$ by the following rule
$$\wt{\s}_w^t(x)=\left[ K(x,t)\wt{\s}_w(x)\right].
$$
\end{de}

It is clear from the very definition that
$$\langle\wt{\s}_w^t(x)\rangle_{\wt I}=\varphi_w(t).
$$
Assume that $\wt{\s}_w^t(x)=\ds\sum_u\alpha_{u,w}(t)\wt{\s}_{w_0u}(x)$.

\begin{lem}\label{l2.2}
\begin{equation}\label{2.4}
\alpha_{u,w}(t)=\frac{\partial}{\partial t_n}\frac{\partial}{\partial t_w}
{\cal F}(t).
\end{equation}
\end{lem}

{\it Proof.}~~ By definition,
\begin{eqnarray*}
\frac{\partial}{\partial t_u}\frac{\partial}{\partial t_w}{\cal F}(t)
&=&\langle\frac{\partial}{\partial t_u}\frac{\partial}{\partial t_w}
K(x,t)\rangle_{\wt I}=\langle\wt{\s}u(x)\wt{\s}_w(x)K(x,t)
\rangle_{\wt I}\\
&=&\langle\wt{\s}_w^t(x)\wt{\s}_u(x)\rangle_{\wt I}=
\langle\wt{\s}_w^t(x),\wt{\s}_u(x)\rangle_Q=\alpha_{u,w}(t).
\end{eqnarray*}
It follows from (\ref{2.3}) that 
$\alpha_{u,w}(t)=\ds\frac{\partial}{\partial t_u}\varphi_w(t)$.

\qed

\begin{lem}\label{l2.3}  For each $w\in S_n$,
\begin{equation}\label{2.5}
\frac{\partial}{\partial t_w}K(x,t)=\wt{\s}_w^t(x).
\end{equation}
\end{lem}

{\it Proof.}~~ We have,
$$\frac{\partial}{\partial t_w}K(x,t)=\left[\wt{\s}_w(x)K(x,t)\right]=
\wt{\s}_w^t(x).
$$
\qed

\begin{cor}\label{c2.1} For each $u,w\in S_n$, let us consider the operator
$$
\Delta_{u,w}=\frac{\partial}{\partial t_u}\frac{\partial}{\partial t_w}-
\sum_{\tau}\wt c_{w,w_0\tau}^{w_0u}\frac{\partial}{\partial t_{\tau}},
$$
where $\wt c_{uv}^w:=\wt c_{uv}^w(q)$ stands for the structural constants
for quantum Schubert polynomials, see Section~2. Then 
\begin{equation}\label{2.6}\Delta_{u,w}{\cal F}=0.
\end{equation}
\end{cor}

{\it Proof.}~~ It follows from (\ref{2.5}), that
\begin{eqnarray*}
\frac{\partial}{\partial t_w}K(x,t)&=&\left[\wt{\s}_w(x)K(x,t)\right]=
\sum_{\tau}\left[\wt{\s}_w(x)\wt{\s}_{w_0\tau}(x)\varphi_{\tau}(t)\right]\\
&=&\sum_u\left[\sum_{\tau}\varphi_{\tau}(t)\wt c_{w,w_0\tau}^{w_0u}\right]
\wt{\s}_{w_0u}.
\end{eqnarray*}
Hence,
$$\alpha_{u,w}(t)=\frac{\partial}{\partial t_w}\varphi_u(t)=\sum_{\tau}
\varphi_{\tau}(t)\wt c_{w,w_0\tau}^{w_0u}.
$$
Equation (\ref{2.6}) follows from the following relations:
\begin{eqnarray*}
\alpha_{u,w}(t)&=&\frac{\partial}{\partial t_u}\frac{\partial}{\partial 
t_w}{\cal F},\ \ {\rm see~Lemma~\ref{l2.2},~ and}\\
\varphi_w(t)&=&\frac{\partial}{\partial t_w}{\cal F}, \ \ {\rm 
see~Lemma~\ref{l2.1}.}
\end{eqnarray*}
\qed

Now we are going to study the structural constants for multiplication of 
the $t$--deformed quantum Schubert polynomials $\wt{\s}_w^t(x)$. Let us 
define functions $\Lambda_{uw\tau}(t)$ from the decomposition
\begin{equation}\label{2.7}
\left[\wt{\s}_u^t\wt{\s}_w^t\right] =\sum_{\tau}\Lambda_{uw\tau}(t)
\wt{\s}_{w_0\tau}^t:=\wt\s_u^t\circ\wt\s_w^t.
\end{equation}

\begin{lem}\label{l2.4}
$$\Lambda_{uw\tau}(t)=\frac{\partial}{\partial t_u}
\frac{\partial}{\partial t_w}\frac{\partial}{\partial t_{\tau}}{\cal F}(t).
$$
\end{lem}

{\it Proof.}~~ By definition,
\begin{eqnarray*}
\wt{\s}_u^t\wt{\s}_w^t&=&\left[\wt{\s}_u\wt{\s}_wK(x,t)K(x,t)\right]
=\sum_{\alpha}\left[\varphi_{w_0\alpha}\wt{\s}_u\wt{\s}_w\wt{\s}_{\alpha}
K(x,t)\right]\\
&=&\sum_{\tau}\left[\sum_{\alpha}\varphi_{w_0\alpha}\langle\wt{\s}_u
\wt{\s}_w\wt{\s}_{\alpha}\wt{\s}_{\tau}\rangle\wt{\s}_{w_0\tau}
K(x,t)\right]=\sum_{\tau}\left(\sum_{\alpha}\varphi_{w_0\alpha}\langle
\wt{\s}_u\wt{\s}_w\wt{\s}_{\alpha}\wt{\s}_{\tau}\rangle\right)
\wt{\s}_{w_0\tau}^t.
\end{eqnarray*}
Hence,
$$\Lambda_{uw\tau}(t)=\sum_{\alpha}\varphi_{w_0\alpha}(t)\langle\wt{\s}_u
\wt{\s}_w\wt{\s}_{\alpha}\wt{\s}_{\tau}\rangle ,\ \ {\rm where}
$$
$$\langle\wt{\s}_u\wt{\s}_w\wt{\s}_{\alpha}\wt{\s}_{\tau}\rangle
=\langle\wt{\s}_u\wt{\s}_w\wt{\s}_{\alpha}\wt{\s}_{\tau}\rangle_{\wt I}.
$$
But by the very construction,
\begin{eqnarray*}
\frac{\partial}{\partial t_u}\frac{\partial}{\partial t_w}
\frac{\partial}{\partial t_{\tau}}{\cal F}&=&\langle\wt{\s}_u\wt{\s}_w
\wt{\s}_{\tau}K(x,t)\rangle_{\wt I}\\
&=&\sum_{\alpha}\varphi_{w_0\alpha}\langle\wt{\s}_u\wt{\s}_w\wt
{\s}_{\tau}\wt{\s}_{\alpha}\rangle_{\wt I}=\Lambda_{uw\tau}(t).
\end{eqnarray*}
\qed

\begin{cor}\label{c2.2} The Gromov--Witten potential ${\cal 
F}(t):=\langle K(x,y)\rangle_{\wt I_n}$ satisfies the 
WDVV--equation (\ref{0.2}).
\end{cor}

{\it Proof.}~~ Follows from associativity of multiplication of 
$t$--deformed Schubert polynomials.
\qed

\begin{cor}\label{c2.3}
$$\langle\wt{\s}_u\wt{\s}_w\wt{\s}_{\tau}\rangle =
\frac{\partial^3}{\partial t_u\partial t_w\partial t_{\tau}}
{\cal F}(t)|_{t=0}.
$$
\end{cor}

{\it Proof.}~~ If $t_w=0$ for all permutations $w\in S_n$, then $K(x,0)=1$,
$\wt{\s}_w^t(x)|_{t=0}=\wt{\s}_w(x)$ and $\Lambda_{uw\tau}(0)=\langle
\wt{\s}_u\wt{\s}_v\wt{\s}_{\tau}\rangle$.

\qed

\begin{lem}\label{l2.5} 
$$\frac{\partial}{\partial t_u}\wt{\s}_w^t(x)=\sum_{\tau}\wt c_{uw}^{\tau}
\wt{\s}_{\tau}^t(x).
$$
\end{lem}
\begin{cor}\label{c2.4} The kernel $K(x,t)$ satisfies the following 
system of differential equations
$$\left(\frac{\partial}{\partial t_u}\frac{\partial}{\partial t_w}-
\sum_{\tau}\wt c_{uw}^{\tau}\frac{\partial}{\partial t_{\tau}}\right)
K(x,t)=0.
$$
\end{cor}

The $t$--deformed quantum Schubert polynomials $\wt\s_w^t$ do not 
orthogonal with respect to the quantum pairing $\langle ,\rangle_Q$ 
any more, but with respect to a
new scalar product $\langle ~,~\rangle_t$ on the quantum cohomology ring:
$$\langle f,g\rangle_t=\langle fgK^{-2}(x,t)\rangle_{\wt I}
$$
they do orthogonal:
\begin{lem}\label{l2.6} The $t$--deformed quantum Schubert polynomials 
$\wt\s_w^t$ are orthogonal with respect to the pairing 
$\langle ~,~\rangle_t$.
\end{lem}

{\it Proof.}~~ It is clear that
$$\langle\wt{\s}_u^t,\wt{\s}_w^t\rangle_t=
\langle\wt{\s}_u^t,\wt{\s}_w^t\rangle_Q=\left\{\begin{array}{rl} 1, & {\rm 
if} \ \ v=w_0w,\\ 0, & {\rm otherwise.}\end{array}\right.
$$
\qed

In the last part of this Section we are going to study the action of 
operators $\ds\frac{\partial}{\partial t_u}$ on the $t$--deformed 
quantum Schubert polynomials $\wt\s_w^t$ (toy analog of Dubrovin's 
connection).
\begin{lem}\label{l2.7} 
$$\frac{\partial}{\partial t_u}\wt{\s}_w^t=\wt{\s}_u\circ\wt{\s}_w^t.
$$
\end{lem}

{\it Proof.}~~ It is clear that
$$\wt{\s}_u\circ\wt{\s}_w^t=\left[\wt{\s}_u\wt{\s}_wK(x,t)\right]=
\sum_{\tau}\Lambda_{uw\tau}(t)\wt{\s}_{w_0\tau}.
$$
Similarly, $\ds\frac{\partial}{\partial t_u}\wt{\s}_w^t=\left[\wt{\s}_u
\wt{\s}_wK(x,t)\right]$.

\qed

One can define a new multiplication $*$ on the ring 
$Q_tH^*(Fl_n)={\rm QH}^*(Fl_n,\Z )\otimes\Q[[t_w]]$:
$$\wt{\s}_u*\wt{\s}_w=\sum_{\tau}\Lambda_{uw\tau}(t)\wt{\s}_{w_0\tau}.
$$

\begin{lem}\label{l2.8} For each $u\in S_n$ let us define an operator 
$\nabla_u$: $Q_tH^*(Fl_n)\to Q_tH^*(Fl_n)$ by the following rule 
$\nabla_u\wt{\s}_w(x)=\ds\frac{\partial}{\partial t_u}\wt{\s}_w^t(x)$. 
Then
$$\nabla_u\wt{\s}_w=\wt{\s}_u*\wt{\s}_w.
$$
\end{lem}

{\it Proof.}~~
$$\nabla_u\wt{\s}_w=\frac{\partial}{\partial t_u}\wt{\s}_w^t=
\sum_{\tau}\Lambda_{uw\tau}(t)\wt{\s}_{w_0\tau}=\wt{\s}_u*\wt{\s}_w.
$$
\qed

\section{Lax pair}
\label{sec:lp}
\neweq

In this Section we construct the Lax pair related to yet another 
deformation $X_w^t$ of the quantum Schubert polynomials.

Define a new scalar product
$$\langle f,g\rangle_t=\langle fgK(x,t)\rangle_{\wt I_n}.
$$

Let $X_w^t(x)$ be Gram--Schmidt's orthogonalization of the 
lexicographically ordered mo\-nomials $x^I$, $I\subset\delta_n$, with 
respect to the pairing $\langle ~,~\rangle_t$, then

\begin{lem}\label{l2.9}
$$\frac{\partial}{\partial t_w}\langle f,g\rangle_t=
\langle\frac{\partial}{\partial t_w}f,g\rangle_t+
\langle f,\frac{\partial}{\partial t_w}g\rangle_t+
\langle\wt{\s}_wf,g\rangle_t.
$$
\end{lem}

Let us define 
\begin{eqnarray}
\wt{\s}_wX_v^t&=\sum_u\varphi^u_{vw}X_u^t \label{3.8}\\
\frac{\partial}{\partial t_w}X_v^t&=\sum_u\psi^u_{vw}X_u^t ,\label{3.9}
\end{eqnarray}
and let us introduce the following matrices

$\bullet$ $L_w=\left((L_w)_{uv}\right)$, where $(L_w)_{uv}=\varphi^v_{uw}(t)$,

$\bullet$ $M_w=\left((M_w)_{uv}\right)$, where 
$(M_w)_{uv}=\psi^v_{uw}(t)$.\\
Then we can rewrite (\ref{3.8}) and (\ref{3.9}) as follows
\begin{eqnarray*}
\wt{\s}_w\cdot X&=&L_w\cdot X,\\
\frac{\partial}{\partial t_w}\cdot X&=&M_w\cdot X,
\end{eqnarray*}
where $X=(X_w^t, w\in W)^t$ is a vector of length $n!$.

\begin{lem}\label{l2.10} Let $u,v\in S_n$, then
$$\frac{\partial}{\partial t_u}L_w=[M_u,L_w]=M_uL_w-L_wM_u.
$$
\end{lem}

{\it Proof.}~~ Let us compute the following expression 
$\ds\frac{\partial}{\partial t_u}\wt{\s}_w\cdot X_v^t$ in two ways, using 
the fact that operators $\ds\frac{\partial}{\partial t_u}$ and $\wt{\s}_w$ 
are commute.

\qed

Let us define ${\cal F}_u(t):=\langle X_u^t,X_u^t\rangle\ne 0$.

\begin{lem}\label{l2.11}
$$\varphi^v_{uw}\cdot{\cal F}_v=\varphi^u_{vw}\cdot{\cal F}_u.
$$
\end{lem}

In other words, if we define $\F ={\rm diag}({\cal F}_v, v\in W)$ then 
matrix $\wt L_w:=L_w{\cal F}$ is symmetric and the matrix $L_w$ is a 
symmetrizable. 

Let us define $\wt M_w=M_w\F$, then

\begin{lem}\label{l2.12}
$$\wt M_w+\wt M_w'+\wt L_w=\frac{\partial}{\partial t_w}\F .
$$
\end{lem}

{\it Proof.}~~ Let us consider $\ds\frac{\partial}{\partial t_w}\langle 
X_u^t,X_v^t\rangle_t=\frac{\partial}{\partial t_w}{\cal 
F}_u\cdot\delta_{uv}$, and apply Lemma~\ref{l2.9}.

\qed

\end{document}